\documentclass{lms}
\usepackage{amssymb,amsmath}
\usepackage{epsfig}
\usepackage{xcolor}
\usepackage{systeme}
\def\RR{\mathbb{R}}
\def\CC{\mathbb{C}}
\def\PP{\mathbb{P}}
\newtheorem{theorem}{Theorem}[section] 
\newtheorem{lemma}[theorem]{Lemma}     

\newnumbered{assertion}{Assertion}    
\newnumbered{conjecture}{Conjecture}  
\newnumbered{definition}{Definition}
\newnumbered{hypothesis}{Hypothesis}
\newnumbered{remark}{Remark}
\newnumbered{note}{Note}
\newnumbered{observation}{Observation}
\newnumbered{problem}{Problem}
\newnumbered{question}{Question}
\newnumbered{algorithm}{Algorithm}
\newnumbered{example}{Example}
\newunnumbered{notation}{Notation} 
\begin{document}
\begin{center}\Large
  {Hyperbolic trigonometric functions as approximation kernels \\ and their properties II: Wavelets}
\end{center}
\bigskip
\begin{center}
  {M. Buhmann, J. J\'odar and M. Rodr\'iguez}
\end{center}
\bigskip
\begin{abstract}
  In a previous paper we have introduced a new class of  radial basis functions
  that are powerful means to approximate functions by quasi-interpolation. In
  this article we extend the results to create new ways of approximating
  functions by prewavelets that are constructed from spaced spanned of the
  new hyperbolic radial basis functions. They consist of highly localised
  time-frequency decompositions that are suitable for analysis and filtering. The construction is
  sufficiently general to apply for large classes of other radial
  basis functions too -- such as multiquadrics and their generalisations and thin-plate splines --, as well as, for example, polynomial splines.
\end{abstract}
\section{Introduction}
There is a plethora of schemes to approximate functions that are
smooth or, for example, represent square-integrable signals by
expansions, interpolation \cite{RBF}, (least-squares) smoothing,
quasi-interpolation \cite{QI}  and wavelet expansions.

These approximations may be based on linear function spaces generated
by polynomials, exponential, piecewise polynomials or radial basis
functions, for example. The latter are particularly interesting if the
nodes on which the data are based are not necessarily equally spaced.

A selection of approximation methods currently available that are closely
related to our work here are
\begin{enumerate}
\item splines (we shall have examples of our results using splines),
\item radial basis functions which are fundamental to this and our previous
  work on hyperbolic trigonometric functions,
  \item wavelet and prewavelet expansions for square-integrable approximands,
\item interpolants and quasi-interpolation as means to generate the approximation of functions and data.
  \end{enumerate}
The approximations make use of (possibly unknown) functions possible for 
\begin{itemize}
\item fast computations, 
evaluations, 
\item compression, data storage, image registration,
\item applications in science, medicine and engineering, where data may come from physical processes
(samples of physical quantities or simulations),  
\item and interpolations (predictions).
\end{itemize}
The authors are keenly interested in the two approaches quasi-interpolation and
wavelets \cite{JJMM} where in  this
article we focus on the wavelet or, in a slightly more general case,
prewavelet part. An earlier article
\cite{HYP} concentrated on the analysis of the Fourier transforms of
freshly introduced new classes of radial functions and their use for
quasi-interpolation. These new classes use hyperbolic tan functions as
radial basis functions which are particularly suitable due to their
local behaviour and their smoothness. Unlike the ubiquitous
multiquadrics or Gau\ss kernel, for example, their smoothness at the
origin and elsewhere is highly adaptable and can be set to different
degrees by judicious choice of parameters.

All constructions are based on creating linear approximation spaces as
spans of shifts of either the radial basis functions themselves or the
quasi-Lagrange functions created thereby. We now wish to construct
prewavelets and linear spaces generated by them in order to decompose
signals of finite energy efficiently and stably.

The detailed knowledge of the properties of Fourier transforms and
quasi-interpolation from \cite{HYP} opens the door to the introduction
of an additional  means to approximate functions and signals:
prewavelets and wavelets \cite{MDBCAM}. Their fundamental idea is to decompose
finite-energy signals in double series of time and frequency, and they
rely on a ``useful'' linear function space and the existence of
square-integrable bases therein. The latter should also be able to perform
basic types of uniform approximation for instance by quasi-interpolation to avail
oneself of wavelet expansions and approximations of sufficient
accuracy.

The ``useful'' spaces (we have conditions for suitable properties
below) are the sets of spaces of square-integrable functions $V_j$,
$j\in\mathbb Z$ and the linear spaces created by the wavelets or
prewavelets are their orthogonal complements $W_j$,
$j\in\mathbb Z$, such that $V_{j+1}=V_j\oplus W_j $ with $W_j\perp
V_j$ for all $j$.

Therefore, in this paper we use the mentioned function spaces spanned
by old and new examples of radial basis functions to construct the
bi-orthogonal basis functions that are required for the development of
wavelets. We include a list of worked-out examples in simple cases
such as cubic B-splines to demonstrate how the methods work and the
procedure to constructing new expansions from the ubiquitous thin-plate
spline and multiquadrics basis functions and of course our new
hyperbolic RBFs from \cite{HYP}.

For intermediate constructions that eventually form the prewavelets, we use
both quasi-interpolation and cardinal interpolation. Quasi-interpolation with radial basis functions (RBFs) is a powerful
technique widely used in numerical analysis and approximation
theory. It provides a flexible and efficient approach for constructing
approximate functions based on given data points without explicitly
performing interpolation. Applications of quasi-interpolation with
RBFs include surface reconstruction, data fitting, image processing,
and solving partial differential equations. It finds utility in
various fields, such as computer graphics, computational physics,
geostatistics and machine learning.

In contrast to this, cardinal interpolation requires to find cardinal
Lagrange functions that interpolate Dirac functions: they have to be one at the
origin and vanish at all other points on an integer grid. They come from the
same linear spaces as the quasi-interpolants but usually require infinite
linear combinations of shifts of radial basis functions.

The article is organised as follows: we begin with some notation and
preliminaries. Then, in the third section, we describe the fundamentals of the construction of
the prewavelets. Finally, there will be an extensive set of examples
which illustrates the method and provides results that can already be
used in applications.

\section{Notations and preliminaries}
There are many ways to introduce the Fourier transform  of a function
$f$ with different scalings. The function must be square-integrable or
absolutely integrable so that its Fourier transform is
square-integrable too  or bounded or -- depending on the integrand's
smoothness -- decaying. The spaces of square-integrable functions
will be denoted by $L^2(\mathbb{R}^n)$, those of square-summable coefficient
sequences by $\ell^2$. We choose 
\[
\hat{f}(\omega)=\int_{\mathbb{R}^n} f(x) e^{-i \omega \cdot x} dx,\quad \omega \in\mathbb{R}^n,
\]
in Jones'  sense (see \cite{JON}). In particular,
let $\cal S$ be the space of infinitely smooth functions
$\gamma:\RR^n\to\RR$  that are compactly supported.
Then the generalised Fourier transform  is that function
$\hat\sigma:{\RR^n}'\to\RR$, $\; {\RR^n}'$ being the origin
punctured space,
of any locally integrable $\sigma:\RR^n\to\RR$ which satisfies
$ \int_{\RR^n}\sigma\hat\gamma= \int_{\RR^n}\hat\sigma\gamma$
for all $\gamma\in\cal S$ which are orthogonal to
$\PP^{k-1}_n$. This denotes the space of polynomials in $n$ variables of total degree $<k$. The hat, when  applied to the test functions whose norm
$\|\gamma\|_1$ is always finite, denotes the classical Fourier transform.
Here $k$ is an integer and the smallest such integer for the $\sigma$
is its infinite discontinuity order.  We shall have examples below.

Our results depend on the
construction of suitable, decaying (square-integrable) basis functions
for our radial basis function spaces. They come along in the shape of
the so-called quasi-Lagrange functions that were originally used for
quasi-interpolation. Their existence and their properties --
especially localness and polynomial reproduction -- are
well-understood by the following theorem.

The following theorem tells us the properties that are sufficient to be satisfied
for a function in order to be a quasi-Lagrange function with certain
polynomial
reproduction properties (at a minimum, to have shifts that form a
partition of unity), see \cite{QI}
or \cite{JJMM} or \cite{MBFD},
for example. The conditions are named after Strang and Fix \cite{SF} who used
them before for the purpose of finite element approximations. In the following form they are generalised for non-compactly supported quasi-Lagrange functions. This requires a substantial change in the conditions for the desired polynomial
reproduction \cite{RBF,QI,MBFD}. We use
the notation $\mathbb{P}_m$ for all polynomials in $n$ variables of
total degree up to and including $m$. The notation $|\cdot|$ is for the
$\ell^1([1\ldots n])$ norm on integer indices.
\begin{theorem}\label{SF}[Strang and Fix conditions]
Let $\psi:\mathbb{R}^n \rightarrow
\mathbb{R}$ be a continuous function such that
\begin{enumerate}
\item there exists a positive real $\ell$ such that for some nonnegative integer $m$, when
  $\|x\|\rightarrow \infty$, $|\psi(x)|={O}(\Vert
  x\Vert^{-n-m-\ell})$, which immediately implies $m$-fold
  differentiability of the Fourier transform,
\item $D^{\alpha} \hat{\psi}(0)=0$, $\forall \alpha \in \mathbb{Z}_+^n$, $1\leq \vert \alpha\vert\leq m$, and $\hat{\psi}(0)=1$,
\item $ D^{\alpha}\hat{\psi}(2\pi j)=0,\ \forall j \in \mathbb{Z}^n\setminus\{0\}$ and  $\forall \alpha \in \mathbb{Z}_+^n$ with $ \vert\alpha\vert \leq m$. 
\end{enumerate}

Then the quasi-interpolant
\begin{equation}
Q_hf(x)=\sum_{j\in\mathbb{Z}^n}f(hj)\psi(x/h-j),\qquad x\in\mathbb{R}^n,
\end{equation}
is well-defined and exact on $\mathbb{P}_m$. The approximation error
for $0<h<1$ can be estimated by
\[\Vert Q_hf-f\Vert_{\infty}=\begin{cases}{O}(h^{m+\ell}), & \text{when } 0<\ell<1,\\
{O}(h^{m+1}\log (1/h)), & \text{when } \ell=1,\\
{O}(h^{m+1}), & \text{when } \ell>1,\\
\end{cases}\]
for $h\rightarrow 0$ and a bounded function $f\in C^{m+1}(\mathbb{R}^n)$ with bounded derivatives.
\end{theorem}
The exactness of the quasi-interpolant on the space ${\mathbb P}_m$ of
course means that whenever $f\in {\mathbb P}_m$, then $Qf\equiv f$.
In our applications (see \cite{QI} or \cite{MBFD}), the
quasi-Lagrange functions are formed as linear combinations of radial
basis functions $\phi: [0,\infty) \rightarrow \mathbb{R}$
\[
\psi(x)=\sum_{\ell\in N}\mu_\ell\phi(|x-\ell |),\qquad x\in{\mathbb
  R},\]
with either a finite $N$ (\cite{MBFD}) or an infinite set $N$
(\cite{JJMM}). This depends on the parity of the singularity of the
radial basis function's generalised Fourier transform at the origin:
if for a nonzero $A$
\[
\hat\phi(r)=A r^{-\mu} + O(r^{-\mu+\epsilon}),\qquad r\to0,
\]
for an even integer $\mu$, the index set $N$ may be finite, otherwise
it will have to be an infinite set (usually, the whole set of
integers). The analysis in the quoted papers then shows that under
suitable extra conditions in the various cases, the largest integer
$m<\mu$ will be admissible.

Important applications and examples of the use of the result above are odd powers
of $r$ as radial basis functions, multiquadrics and its
generalisations \cite{ORT1,ORT}, thin-plate splines which are according to (Buhmann and Ortmann, 2025) \cite{ORT}
\begin{equation}\label{OBno2}
\phi(r)=\Bigl(r^{2\beta}+b^{2\beta}\Bigr)^\gamma\log \Bigl(r^{2\beta}+b^{2\beta}\Bigr)\end{equation}
in the most general form,
and their
brethren and our novel radial basis functions from \cite{HYP}.

All of those examples depend on the order of the singularities of the radial
basis function's generalised Fourier transform at the origin. For
details see \cite{QI} and \cite{JON} and \cite{ORT1,ORT}.

The key features of the approximation with quasi-interpolation and 
their approximation orders (approximation powers) stem from suitable conditions. 
These conditions and indeed the obtainable approximation
powers depend essentially on the underlying space dimensions and on the choice of
the radial basis function. For the latter, the choices of multiquadrics in the
currently most general form of Buhmann and Ortmann (2024) \cite{ORT1}
\begin{equation}\label{OBno1}\phi(r)=\Bigl( r^{2\beta}+b^{2\beta}\Bigr)^{\gamma/2},
\end{equation}
inverse 
multiquadrics (this is (\ref{OBno1}) for $\beta=1$, $\gamma=-1$) and inverse quadratics which is (\ref{OBno1}) for $\beta=1$ and $\gamma=-2$, odd powers of $r$ and even powers of $r$ times 
a logarithm (\ref{OBno2}) are some of the most interesting choices.

The expansions of these radial functions and their distributional Fourier transforms being well-understood
we can use classical means such as Poisson's summation formula to regard the
Fourier transforms of the linear combinations of the shifts of then radialised
functions in order to form and satisfy the aforementioned conditions on the linear 
combination's  coefficients.   

We should note already now that some RB-function-kernels do not admit our constructions
of quasi-interpolants, such as the following or in the following situations: if the radial functions
are absolutely integrable or even compactly supported or exponentially
decaying, their Fourier transforms exist in the classical way and are
certainly continuous about the origin and sometimes even
analytic. When the Fourier transforms have no infinite discontinuity, and the
original radial functions have a nonzero integral,
this Fourier  transform has a well-defined nonzero value at the origin
and this admits not linear combinations of translates that would  then
be used to form the pre-conditioned-interpolants.

The radial basis functions mentioned before the previous paragraph admit the application of Theorem \ref{SF} in one dimension for 
\begin{enumerate}[(i)]
\item $\mu=2,m=1$ for multiquadrics $\phi(r)=\sqrt{r^2+b^2}$, i.e., $\beta=\gamma=1$ in (\ref{OBno1}),
\item $m=0$ for inverse multiquadrics $\phi(r)=1/\sqrt{r^2+b^2}$, i.e., $\beta=1$, $\gamma=-1$ in (\ref{OBno1}) --
  this is a special case, where a logarithmic singularity comes up in its generalised Fourier transform 
  (see \cite{RBF}),
\item $\mu=3, m=2$ for thin-plate splines \cite{JJMM} $\phi(r)=r^2\log r$, i.e., $\beta=\gamma=1$ in (\ref{OBno2})
  \item $m=1$ for the new hyperbolic tangents radial functions
    $\phi(r)=r\tanh r$, \cite{HYP}, and
    \item $\mu=3, m=2$ for the generalised multiquadrics
      $\phi(r)=\sqrt{r^4+c^4}$ \cite{ORT1}. Those are the cases of the radial
      basis functions (\ref{OBno1}) of Ortmann and Buhmann with $\beta=2$ and $\gamma=1$.
  
  \end{enumerate}

\section{General construction}
There are many equivalent ways for defining wavelets and
prewavelets. They are usually based on the definition of a
multiresolution analysis. We shall also use this form of creating
wavelet spaces and wavelet decompositions with our new radial basis functions.
The basic idea is to decompose $L^2({\mathbb R})$ into nested
sequences of spaces that become dense in the space and have an
intersection of zero only. To wit,
\begin{definition} A multiresolution analysis (MRA) of $L^2(\mathbb{R})$ consists of a nested sequence of closed
subspaces $\{V_j\}_{j\in \mathbb{Z}}$ of $L^2(\mathbb{R})$  and a function $\psi \in V_0$, called the scaling function,
such that
\begin{enumerate} 
\item $V_j \subset V_{j+1}$ for  all $j\in \mathbb{Z}$,
\item $ \displaystyle \overline{\bigcup_{j\in \mathbb{Z}} V_{j}}= L^2(\mathbb{R})$,
\item $\displaystyle \bigcap_{j\in \mathbb{Z}} V_{j}= \{0\}$,
\item  $f \in  V_0$ if and only if $f(2^{j}\cdot) \in  V_j$.
\item  The collection
$\{\psi(x-k)\}_{k\in \mathbb{Z}}$ is a Riesz basis for $V_0$, that is
  there are positive and finite constants $A$ and $B$ such that
  \[
A\sum_{\ell=-\infty}^\infty
  |c_\ell|^2\leq\Big\|\sum_{\ell=-\infty}^\infty
  c_\ell\psi(\cdot-\ell)\Big\|^2_2\leq B\sum_{\ell=-\infty}^\infty
  |c_\ell|^2.
  \]
\end{enumerate} 
\end{definition}
We recall that for example B-splines form a Riesz basis for $V_0$ the
linear space of all square-integrable splines of some order.

An orthonormal basis is a Riesz basis for $V_0$ with all constants $A=B=1$.
Any MRA is completely determined by the scaling function $\psi:\mathbb{R}\to\mathbb{R}$.

It is well known that given a MRA, it is straightforward that an
orthonormal wavelet basis of $L^2(\mathbb{R})$ can be constructed. If
the orthogonality between different translates only of the wavelets is
given up, the resulting functions are called prewavelets. We shall use
prewavelets only in this article. It is also known that
\begin{lemma}
Given $f\in L^2(\mathbb{R})$, the system $\{T_jf(x):=f(x-j)\}_{j\in \mathbb{Z}}$ is an orthonormal system if and only if
\[ \displaystyle \sum_{j \in \mathbb{Z}} |\hat{f} (\omega+ 2\pi j
)|^2\equiv1,\qquad \omega\in{\mathbb T}=[-\pi,\pi].\]
\end{lemma}
The proof of this is a straightforward consequence of Poisson's summation formula
\cite{STEINWEISS}.

The idea now is to decompose each $V_{j+1}$ for all integers $j$ into
an orthogonal and direct sum of
\begin{equation}\label{decomp1} V_{j+1}=V_j\oplus W_j\end{equation}
where $W_j\subset V_{j+1}$ and $W_j\perp V_j$. Of course, given the
concept of a multiresolution analysis, it is sufficient to carry this
out for $j=1$ and all other decompositions follows automatically by
scaling. Eventually, all of $L^2(\mathbb{R})$ can be decomposed into
\begin{equation} \label{decomp2} \cdots \oplus W_{-2}\oplus W_{-1}\oplus W_{0}\oplus W_{1}\oplus
W_{2}\cdots.\end{equation}
It is important to point out that (\ref{decomp1}) and (\ref{decomp2})
can also be achieved when condition (iv) of our multiresolution
analysis does not hold, i.e., when we have a sequence of subsets $V_j$
of square-integrable functions for $j\in\mathbb Z$ that are dense in
$L^2$ and whose intersection is $\{0\}$, but then we have to establish
(\ref{decomp1}) for each index separately--and not just for one (and then scaling).

There exists a function $\zeta \in W_0$ which is in the orthogonal
complement of  $V_0$ in $V_1$. And this function which we call the
prewavelet generates $W_0$.
A function $\zeta\in L^2 (\mathbb{R})$   is called a wavelet if $\{\zeta_{j,k} \; :\;  j,k\in \mathbb{Z} \}$ is an orthonormal
basis in $L^2(\mathbb{R})$, where $\zeta_{j, k}(x)= 2^j\zeta(2^j x
-k)$ and $x\in\mathbb R$. A function $\zeta\in L^2 (\mathbb{R})$   is called a prewavelet if $\{\zeta_{j,k} \; :\;  j,k\in \mathbb{Z} \}$ is an orthonormal
basis in $L^2(\mathbb{R})$ with orthogonality conditions with respect
to $j$ but no orthogonality conditions with
respect to $k$. We will only deal with prewavelets in this text.

Our goal is  to define a function, the prewavelet,  $\zeta$,  and  detail  the scheme  of the multiresolution analysis. More precisely, we  set 
$\zeta(x):= \psi(2 x)- \chi(2 x)$ and $x\in\mathbb R$,
where  $\psi$ is the quasi-Lagrange function based on a suitable function $\phi$ and  $\chi$ is the cardinal function based on the same function (see \cite{RBF} for example, for more details on the cardinal function).
In this way, as a start,  $\zeta \in W_0$ and $\psi \in V_0$.

Now, let us assume that we have a function $\psi$ that satisfies all
conditions of the Strang and Fix theorem for positive $m$ and has the
form
\[
\psi(x)=\sum_{\ell\in N}\mu_\ell\phi(|x-\ell |),\qquad x\in{\mathbb R},
\]
for a radial basis function with positive generalised Fourier
transform such that also the aforementioned Lagrange function
\[
\chi(x)=\sum_{\ell\in \mathbb Z}\lambda_\ell\phi(| x-\ell |),\qquad x\in{\mathbb R},
\]
exists. Suppose also the latter is square-integrable -- that the
former is follows from the Strang and Fix theorem. Defining
\begin{equation}\label{vzero} V_0:={\rm span}\{\psi(\cdot-\ell)\mid\ell\in{\mathbb Z}\}\end{equation}
and
\begin{equation}\label{vone} V_1:={\rm span}\{\psi(2\cdot-\ell)\mid\ell\in{\mathbb Z}\},\end{equation}
where in both cases all linear combinations to generate the span are
formed with square-summable coefficients, we need to show that
\begin{equation}\label{wzero} W_0:={\rm span}\{\zeta(\cdot-\ell)\mid\ell\in{\mathbb Z}\}\end{equation}
satisfies $W_0\perp V_0$ and $V_1=V_0\oplus W_0$. This can be done by
solving a certain bi-infinite system of linear equations which we will
do next. We point out that while the $V_k$ become dense in $L^2(\RR)$ because of Theorem \ref{SF}, they are not necessarily nested.
\subsection{Obtaining the $a_{\cdot, i}$ coefficients}
In order to show the required properties of $W_0$ and the prewavelet
we claim that the following two identities hold. Namely,
\begin{equation}\label{struc}
\begin{aligned}
\psi(2x)&=\sum_{k\in \mathbb{Z}} b_{1,k}\psi(x-k)+ \sum_{k\in \mathbb{Z}} a_{1,k}\zeta(x-k) \\
\psi(2x-1)&=\sum_{k\in \mathbb{Z}} b_{2,k}\psi(x-k)+ \sum_{k\in \mathbb{Z}} a_{2,k}\zeta(x-k).
\end{aligned}
\end{equation}
If that is so and if the four sequences of coefficients therein are
all absolutely summable, it follows by H\"older's inequality that the
space $W_0$ generated by all even and odd shifts of
$\psi(2\cdot-\ell)$, $\ell$ all integers, even or odd, satisfies the required
properties. This is because then not only the functions on the
left-hand sides of (\ref{struc}) can be re-written as expansions from
$W_0$ and $V_0$ but also sums of them with square-summable
coefficients. And therefore all functions from the linear space $V_1$
will be included.
\begin{theorem} Let $\psi$ satisfy the Strang and Fix conditions for $n=1$ and $m>0$. Then there is a square-integrable function $\zeta$ in $V_1$ and there are absolutely summable coefficients such that the equations
  (\ref{decomp1}) and (\ref{decomp2}) hold and so that the spaces
 {\rm (\ref{vzero})}, {\rm (\ref{vone})} and {\rm (\ref{wzero})} satisfy $W_0\perp V_0$ and $V_1=V_0\oplus W_0$.
  \end{theorem}

We therefore need to compute those coefficients and
study their summability properties, and we begin with the $a$s.

We start working on the first equation. This is a bi-infinite linear
equation whose matrix corresponds to a bi-infinite operator. The
matrix is a Toeplitz matrix. The Toeplitz matrix is invertible, i.e.,
the operator between $\ell^2({\mathbb Z})\to \ell^2({\mathbb Z})$ is
invertible if the symbol of the Toeplitz matrix forms an absolutely
convergent series that is positive. So, taking inner products with $\zeta(x-\ell)$ on the whole expression we have
\[
\left(\psi(2x), \zeta(x-\ell)\right) =\left(\sum_{k\in \mathbb{Z}} b_{1,k}\psi(x-k), \zeta(x-\ell)\right)+ \left(\sum_{k\in \mathbb{Z}} a_{1,k}\zeta(x-k), \zeta(x-\ell)\right)
\]
for any real $x$. 
Now, as the first term of the right hand side is 0 because of the orthogonality between $\zeta$ and the basis functions $\psi(\cdot - j), j\in \mathbb{Z}$, we have
\begin{equation}\label{a1}
\left(\psi(2x), \zeta(x-\ell)\right) =\sum_{k\in \mathbb{Z}} a_{1,k} \left( \zeta(x-k), \zeta(x-\ell)\right). 
\end{equation}
Our goal is to get the coefficients $a_{1,k}.$ For getting those
coefficients, we will apply the Wiener's lemma \cite{Wiener} and we
must work on the aforementioned symbol, that is, show that it is an
absolutely convergent Fourier series that does not vanish. Then the
coefficients of the inverse of the Toeplitz matrix are the Fourier
coefficients of the reciprocal of the symbol.

The symbol is defined in our case by
\[
\sigma(\vartheta) =\sum_{\ell\in \mathbb{Z}} e^{-i \vartheta \ell }\left( \zeta(x-\ell), \zeta(x)\right)=
\left(  \sum_{\ell\in \mathbb{Z}} e^{-i \vartheta \ell }\zeta(x-\ell), \zeta(x)\right). 
\]
Now, we compute the Fourier transform,  as a function of $y,$ of
$e^{-i \vartheta y} \zeta(x-y)=: g(y)$. We need this in order to apply
Poisson's summation formula in order to rephrase the symbol and study
its properties:
\begin{align*}
\hat{g}(\omega)&=\int_{-\infty}^{\infty} e^{-i \omega y} e^{-i \vartheta y}  \zeta(x-y)dy
= \int_{-\infty}^{\infty} e^{-i y (\omega+  \vartheta)} \zeta(x-y)dy\\
&= \int_{-\infty}^{\infty} e^{-i (x-s) (\omega+  \vartheta)} \zeta(s)ds
=e^{-i x (\omega+  \vartheta)} \int_{-\infty}^{\infty} e^{-is (-\omega-  \vartheta)} \zeta(s)ds\\
&=e^{-i x (\omega+  \vartheta)} \hat{ \zeta}(-\omega-  \vartheta).
\end{align*}
The Poisson summation formula (see \cite[pg.11]{QI} for instance) states that
\[\sum_{\ell\in \mathbb{Z}} e^{-i \vartheta \ell }\zeta(x-\ell)= \sum_{m\in \mathbb{Z}} e^{-i x(2\pi m + \vartheta)} \hat{\zeta} (-2\pi m - \vartheta).\]
So this allows a reformulation of the symbol as 
\begin{align*}
\sigma(\vartheta) &=\left(  \sum_{\ell\in \mathbb{Z}} e^{-i \vartheta \ell }\zeta(x-\ell), \zeta(x)\right)
 =\left( \sum_{m\in \mathbb{Z}} e^{-i x(2\pi m + \vartheta)} \hat{\zeta} (-2\pi m - \vartheta), \zeta(x)\right)\\
&= \int_{-\infty}^{\infty} \sum_{m\in \mathbb{Z}}e^{-i x(2\pi m + \vartheta)}   \hat{\zeta} (-2\pi m - \vartheta) \zeta(x)d x
=  \sum_{m\in \mathbb{Z}} \int_{-\infty}^{\infty} \hat{\zeta} (-2\pi m - \vartheta)  e^{-i x(2\pi m + \vartheta)}  \zeta(x)d x\\
&= \sum_{m\in \mathbb{Z}}\hat{\zeta} (-2\pi m - \vartheta) \hat{\zeta} (2\pi m + \vartheta) =  
\sum_{m\in \mathbb{Z}} |\hat{\zeta} (\vartheta +2\pi m ) |^2\geq 0.
\end{align*}
In order to apply Wiener's Lemma we require two conditions: the symbol must be positive and the coefficients must be in $\ell_1.$
This latter condition holds because $\psi$ and $\chi$  decay
sufficiently fast.  If it is additionally true that $\sigma(\vartheta) >0$, then we have that
\[ \dfrac{1}{\sigma(\vartheta)}=   \sum_{j\in \mathbb{Z}} d_j e^{-i j \vartheta  }\] and we can obtain the coefficients $a_{i,1}$ from (\ref{a1}) and they are
absolutely summable.

The remaining problem is what happens if the symbol vanishes. As 
$\zeta(x)= \psi(2 x)- \chi(2 x)$  taking Fourier transforms on both sides, we have that 
\[\hat{\zeta}(\omega)= \dfrac12 \hat{\psi}\left(\frac{\omega}{2}\right)-\dfrac12\hat{ \chi}\left(\frac{\omega}{2}\right),\]
which gives us
\[
\sigma(\vartheta)= \sum_{m\in \mathbb{Z}} |\hat{\zeta} (\vartheta +2\pi m ) |^2
=\frac14 \sum_{m\in \mathbb{Z}} \Big|\hat{\psi}\left(\frac{\omega}{2}+\pi m\right)-\hat{ \chi}\left(\frac{\omega}{2}+\pi m\right)\Big|^2.
\]
Now, if the symbol vanished at one point, $\vartheta_0,$ then
$\hat{\psi} $ and $\hat{\chi} $ would have to satisfy the equation
\begin{equation}\label{van}
\hat{\psi} (\vartheta_0+ \pi m)=  \hat{\chi} (\vartheta_0+ \pi m), \quad \forall m\in  \mathbb{Z}.
\end{equation}
Now, we are considering that $\psi$ and $\chi$ satisfy 
\[ \hat{\psi}(\cdot)= P(\cdot) \hat{\phi}(|\cdot|), \qquad
\hat{\chi}(\cdot)= \dfrac{\hat{\phi}(|\cdot|)}{\displaystyle \sum_{\ell \in
    \mathbb{Z}}  \hat{\phi}(|\cdot+ 2\pi \ell | )}\] for a certain 
function $\phi$ where $P$ is a trigonometric expansion satisfying
$P(0)=P(2\pi)=0$.

So, taking into account (\ref{van}), we obtain that
\[  P(\vartheta_0+ \pi m)= \dfrac{1}{\displaystyle \sum_{\ell\in \mathbb{Z}}  \hat{\phi}(|\vartheta_0+ \pi m+  2\pi \ell |)}, \quad \forall m\in  \mathbb{Z},\]
which means that for $m$ even we would have 
\begin{equation}\label{even}
 P(\vartheta_0+  \pi m)= P(\vartheta_0)= \dfrac{1}{\displaystyle \sum_{\ell \in \mathbb{Z}}  \hat{\phi}( |\vartheta_0+ 2\pi \ell| )},
 \end{equation}
and when $m$ is  odd 
\begin{equation}\label{odd}
 P(\vartheta_0+  \pi m)= P(\vartheta_0+ \pi ) = \dfrac{1}{\displaystyle \sum_{\ell\in \mathbb{Z}}  \hat{\phi}(|\vartheta_0+ \pi + 2\pi \ell| )}.
\end{equation}
Now, if we take 
\[P_a(\vartheta):= P(\vartheta)+ A(1-\cos(\vartheta))^{\alpha}, \qquad A > 0, \alpha \in \mathbb{N},\] 
for certain values of $A$  and $\alpha$ we would have a new function
$\phi_1(\cdot)$ meeting the identical conditions we need (i.e., those of the
Strang Fix conditions) at the $2\pi$-multiples of all integers for this construction, but the symbol of this new $\phi_1$ would be positive.

The same process can be developed for obtaining the $a_{2, i}$ coefficients  of the second expression  in (\ref{struc}).
\subsection{Obtaining the $b_{\cdot, i}$ coefficients}
Next, we need to compute the coefficients $b$. The process is very similar to the one we have  presented  before. We do some computations in order to
demonstrate the argument.
We multiply the  first equation  of (\ref{struc}) by  $\psi(x-\ell)$ and taking inner products we obtain
\[
\left(\psi(2x), \psi(x-\ell)\right) =\left(\sum_{k\in \mathbb{Z}} b_{1,k}\psi(x-k), \psi(x-\ell)\right)+ \left(\sum_{k\in \mathbb{Z}} a_{1,k}\zeta(x-k), \psi(x-\ell)\right). 
\]
Now, due to the orthogonality conditions,  the second term of the right-hand side vanishes and we obtain by orthogonality
\begin{equation}\label{b1}
\left(\psi(2x), \psi(x-\ell)\right) =\sum_{k\in \mathbb{Z}} b_{1,k} \left( \psi(x-k), \psi(x-\ell)\right). 
\end{equation}
For getting the $b_{1,k} $  coefficients, we must compute  the symbol:
\[
\sigma_1(\vartheta) =\sum_{\ell\in \mathbb{Z}} e^{-i \vartheta l }\left( \psi(x-\ell), \psi(x)\right)=
\left(  \sum_{\ell\in \mathbb{Z}} e^{-i \vartheta l }\psi(x-\ell), \psi(x)\right). 
\]
and we reach, in the same way that we did previously, the expression
\begin{align*}
\sigma_1(\vartheta) =
\sum_{m\in \mathbb{Z}} |\hat{\psi} (\vartheta +2\pi m ) |^2\geq 0,
\end{align*}
so the symbol is non-negative.
Again we can apply the Wiener's Lemma taking care that the symbol must be greater than 0.
For similar reasons, if this is not so, we can replace $P$ by
\[P_b(\vartheta):= P(\vartheta)+ B(1-\cos(\vartheta))^{\beta}, \quad B > 0, \beta \in \mathbb{N}.\]
On one side  if we set $\vartheta^{\ast}=\{0, 2 \pi\}$ we  have that
\[P(\vartheta^{\ast}) = \dfrac{1}{\displaystyle\sum_{\ell \in \mathbb{Z}} |\hat{\psi} (2\pi \ell )|}\]
On the other side, we would have that 
\[P_b(\pi)=P(\pi)+ 2^{\beta} B \neq  \dfrac{1}{\displaystyle \sum_{\ell \in \mathbb{Z}} |\hat{\psi} ( 2\pi \ell )|}.\]
So, both conditions cannot be satisfied at the same time and the symbol must always be positive.

\section{Examples}
Let us consider the cutoff function $(x)_+ = \max\{x, 0\}$ which we
shall need to form our first example, namely the cubic B-splines in one
variable and a prewavelet for a
linear space of cubic splines.

In order to form a B-spline basis for this space,  we recall the
definition of a divided difference (DD) first:
\begin{definition}
With respect to one center $t_0$, a divided difference (DD) is defined by convention as
$[t_0]f = f(t_0),$ and  given a positive integer $k$ the divided difference with respect  $k + 1$ centers $t_0, \ldots, t_{k}$ is defined as
\[
[t_0, \ldots, t_{k}]f =\dfrac{ [t_1, \ldots, t_{k}] f- [t_0, \ldots, t_{k-1}]f}{t_{k} -t_0}.
\]
Given a set of  $n$ real distinct points $\{x_0, x_1, \ldots, x_n\}$ the $n$th divided difference of a function $f$ has the following expression
\[
[x_0, x_1, \ldots,x_n]f= \displaystyle  \sum_{k=0}^n f(x_k)\prod_{\substack{\ell=0 \\ \ell\neq k}}^{n}\dfrac{1}{x_k-x_\ell}.
\]
Note that the divided differences remain the same if the evaluation points are permuted.
\end{definition}
We also recall the standard definition of a B-spline:
\begin{definition}
For $x_0 < \cdots < x_p$, the B-spline of order $p$ with knots $x_0, \ldots, x_p$ is defined by
\[
\RR \ni u \quad\to\quad B_p(u | x_0,\ldots,x_p) = [x_0,\ldots,x_p]\left \{ \frac{(\,\cdot\, - u)_+^{p-1}} {(m-1)!}\right\}.
\]
\end{definition}
The B-spline vanishes outside the interval $(x_0,x_p)$ by definition, and it is strictly positive on the interval itself. Moreover,
\[
     \int_\RR B_p(u|x_0,\ldots,x_p)\;\,d u = \frac{1}{m!}.
\]
If
$f :\RR \to \CC$ is $p$-times continuously differentiable, then also
\[
[x_0,\ldots,x_p]\,f = \int_\RR f^{(p)}(u)B_p(u | x_0,\ldots,x_p)\,  du.
\]
In the notation below we shall use equispaced knots and omit the notation for them in the general form of the B-splines, i.e., 
\[B_p(x)=B_p\Bigl(x | -\frac{p}{2}, -\frac{p}{2} +1, \ldots,\frac{p}{2}\Bigr), \qquad x\in\RR,
\] 
and its Fourier transform is given by
\begin{equation}\label{fou}
\hat{B_p}(\omega)= \left(\frac{\sin \Bigl(\frac{\omega}{2}\Bigr)}{\frac{\omega}{2}}\right)^p = \Bigl(\text{sinc}\left(\frac{\omega}{2}\right)\Bigr)^p, \qquad \omega\in{\mathbb R}.
\end{equation}

In this section we show some examples in order to gain intuition for the theoretical results obtained.  In each of the examples, both the cardinal function, $\chi(\cdot)$, and the quasi-Lagrange function, $\psi(\cdot)$, associated a given function radial basis function $\phi$ are computed. From these, the wavelet $\zeta(\cdot)$ is obtained. All computations to obtain both Fourier coefficients and graphs of wavelet functions have been performed using Mathematica software.
\subsection{Wavelet based on the truncated power  of order 3}
\begin{figure}[ht]
        \centering
        \includegraphics[width=0.4\textwidth]{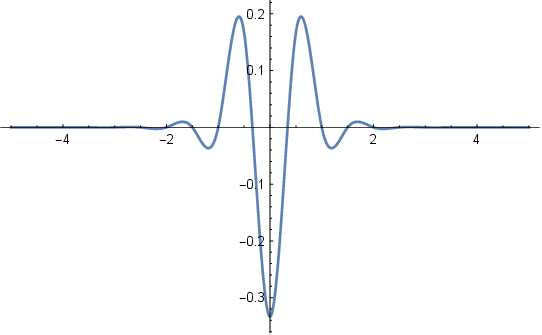}
 \caption{Wavelet based on the truncated power of order 3.}
 \label{alhambra}  
\end{figure}

In this first example we have chosen the function space for approximation and for the MRA as the space of cubic splines; as the scaling function $\phi$ we therefore take the truncated power of order three. 
\begin{description}
\item[Quasi-Lagrange function]
Firstly, we can observe that the cubic B-spline $B_4(x)$ can be a good candidate for being  the quasi-Lagrange function. In order to check  the conditions of Theorem \ref{SF} we can see that (i) is straightforward because of the compact support of the B-splines.  Now, by (\ref{fou}) the Fourier transform of $B_4(x)$ is 
\[\left(\text{sinc}\Bigl(\frac{\omega}{2}\Bigr)\right)^4,\qquad \omega\in{\mathbb R},\]
and it is easily checked  that the conditions (ii)--(iii) of the aforementioned Theorem  hold for $m=1$ and any $\ell$. Therefore we take $\psi:= B_{4}$. 

\item[Cardinal function]
For practical purposes, instead of generating the cardinal function from a linear combination of shifts of  $\phi$, we work with the quasi-Lagrange function $\psi$. There is no problem with this, since $\psi$ is a finite linear combination of shifts of $\phi$. The advantage of doing so is that $\psi$  has compact support, while $\phi$  does not. In this way, it follows from \cite[Sec. 19.3]{POW} that the cardinal function can be expressed in the form 
\[\chi(x)=\sum _{p\in \mathbb{Z}} \lambda_p \psi(x-p),\qquad x\in{\mathbb R},\]
where 
\[\lambda_p= \sqrt{3} \left(\sqrt{3}-2\right)^{\left| p+2\right|},
\quad p\in \mathbb{Z}.\] 
These values come from the resolution of an infinite linear system by imposing the cardinal function to be one at the origin and zero at all other integer. We will get into the details in the last example.

The asymptotic decay of $\chi (x)$ is identified as 
\[
|\chi (x)|={O}((1+|x|)^{-5}),
\]
for large $|x|$ (see, for instance, \cite[Theorem 4.3]{RBF}).
\end{description}
Finally, the wavelet is then given by 
\[  \zeta(x)= \psi(2x)- \chi(2x).\]
The plot of  $ \zeta$ is given in Figure \ref{alhambra}.
\subsection{Wavelet based on the multiquadric function}
For this example, we  consider $\phi$ as the multiquadric function
given by $\phi(r)=-\sqrt{c^2+ r^2}$, with $c>0$. We note that the
standard multiquadric function is provided with a negative sign here
because then the Fourier transform is positive rather than
negative. (This plays no r\^ole of course when interpolation is
carried out.)

We know from \cite{RBF} that the generalised Fourier transform of the multiquadric function is
\[\hat{\phi}(|\omega|)=\frac{2 c}{|\omega|}K_{1}(c  |\omega| ), \qquad \omega\in \mathbb{R},\]
where $K_{1}(\cdot)$ is the modified Bessel function of the second
 kind, also called Macdonald function \cite{AS}.

With this information we can compute the wavelet.
\begin{description}

\item[Quasi-Lagrange function]
The  quasi-Lagrange function is obtained as the  divided difference at
the nodes $-1$, $0$ and $1$ of the multiquadric function (\cite[Chapter 2]{RBF}):
\[\psi(x)=-\frac{1}{2c} \phi(\left| x-1\right| )+\frac{1}{c}\phi(\left| x\right|
)-\frac{1}{2c}\phi(\left| x+1\right| ),  \qquad x\in \mathbb{R}.\]
With this choice, it is satisfied the decay condition (\cite[Proposition 2.3]{RBF})
\[
|\psi (x)| ={O}(|x|^{-3}), \quad |x|\rightarrow \infty.
\]

On the other hand, by computing the Fourier transform of $\psi$, it is obtained that $\hat{\psi}(\omega)=P(\omega) \hat{\phi}(|\omega|)$, where $P$ is the trigonometric polynomial coming from the second divided difference, that is, 
\[
P(\omega)=\frac{1}{c}(1-\cos \omega).
\]
Taking into account that the modified Bessel function $K_1$ satisfies $K_1 (z) \sim z^{-1}$, $z \rightarrow 0_+$ (\cite{AS}), it can be checked that all the conditions in Theorem \ref{SF} are verified for $\psi$ with $m=\ell=1$.

\item[Cardinal function]

To obtain the cardinal function, $\chi$, we need to know the symbol function, $\sigma$, and from it calculate the Fourier coefficients of the reciprocal of the symbol, that is,
\[\chi(x) := \sum _{k\in \mathbb{Z}}  c_k \phi(\left| x-k\right|
),\qquad x\in\RR,\]
where
\[
c_k=\frac{1}{2\pi} \int_{-\pi}^{\pi} \frac{e^{{\rm i}\omega k}}{\sigma (\omega)} d\omega.
\]

 In our case the symbol is given by
\[\sigma (\omega)=\sum _{j\in \mathbb{Z} } \hat{\phi} (\left| \omega+ 2 \pi  j \right|)=2c \sum _{j\in \mathbb{Z}}  \frac{K_1(c \left| \omega+ 2 \pi  j\right| )}{\left| \omega+ 2 \pi  j \right|}, \qquad \omega\in \mathbb{R}.\] 
  
 It should be noticed that $\chi$ also decays cubically at infinity as $\psi$ did (\cite[Theorem 2.7]{RBF}.

To perform the computations, as the tails of the symbol function go rapidly to zero, we have truncated the infinite sum that is involved. Specifically, we have truncated the sum taking $-5\leq j \leq 5$ for the index. From this truncated sum  we have computed a set $\Delta$ of $2^{10}$ Fourier coefficients of the reciprocal of the symbol in order to obtain an approximated version of the cardinal function
\[\chi_t(x) := \sum _{c_k\in \Delta}  c_k \phi(\left| x-k\right|
),\qquad x\in\RR.\]
 We should note that not for all integers $i$  we have $\chi_t (i)= \delta_{0i}$ due to this double truncation (first of the symbol and then in the calculation of Fourier coefficients). However, the approximation $\chi_t (i) \approx \delta_{0i}$  is quite accurate.

\end{description}

\begin{figure}[h] 
        \centering
        \includegraphics[width=0.4\textwidth]{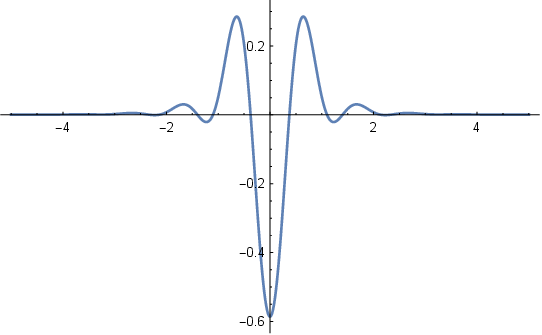}
 \caption{Wavelet based on the multiquadric function.}
 \label{leine} 
\end{figure}

The wavelet is given by 
\[  \zeta(x)= \psi(2x)- \chi(2x), \qquad x\in \mathbb{R}.\] 
and its plot is given in Figure \ref{leine}.
\subsection{Wavelet based on a hyperbolicus function}\label{ex:hyp}
We now take $\phi(r)=- r \tanh r$. The generalised Fourier transform of the function $\phi$ is
given (see the appendix for the computation) by
\begin{equation}\label{HypFT}
\hat{\phi}(\omega)=\frac{ \pi ^2}{2} \coth \left(\frac{\pi
  \omega}{2}\right) \text{csch} \left(\frac{\pi  \omega}{2}\right) =
\frac{ \pi ^2}{2} \frac{\cosh \left(\frac{\pi
    \omega}{2}\right)}{\sinh^2 \left(\frac{\pi  \omega}{2}\right)}, \qquad \omega\in \mathbb{R}.
\end{equation}
    
\begin{figure}[h]
        \centering
        \includegraphics[width=0.4\textwidth]{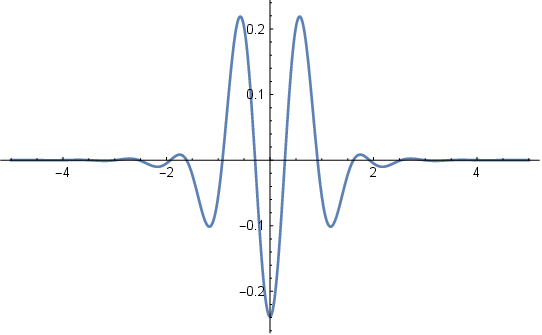}
 \caption{Wavelet based on $\phi(r)=- r \tanh r$.}
 \label{canena}  
\end{figure}

\begin{description}

\item[Quasi-Lagrange function]
This generalised Fourier transform $\hat{\phi}$ has  a singularity of order two  at the origin since
 \[
 \hat{\phi}(|\omega|) \sim \frac{2}{\omega ^2}, \qquad |\omega|\rightarrow 0.
 \]
 
The quasi-Lagrange function will be given by
\[
\psi(x)=\sum_{k \in \mathbb{Z}} \mu_k \phi (|x-k|),
\] 
where the $\mu_k$ are the Fourier coefficients of the trigonometric polynomial $P(\omega)=(1-\cos \omega)$, for instance. With this choice, it is obtained that
\[
|\psi(x)|={O}(|x|^{-3}),
\]
for large $|x|$, and also that 
\[
\hat{\psi}(\omega)=P(\omega) \hat{\phi}(|\omega|)=1+{O}(|\omega|^2)
\]
(see \cite{HYP}). Therefore, the conditions in Theorem \ref{SF} are satisfied with $m=\ell=1$. Concretely, the  quasi-Lagrange function which is obtained is, as in the previous example, the divided difference at the nodes $-1$, $0$ and $1$ of the function $\phi(|\cdot|)$:
\[\psi(x)=-\frac{1}{2} \phi(\left| x-1\right| )+\phi(\left| x\right| )-\frac{1}{2}\phi(\left| x+1\right| ), \qquad x\in \mathbb{R}.\]
\item[Cardinal function]
Taking into account (\ref{HypFT}), the symbol is now given by 
\[
\sigma (\omega)= \frac{ \pi ^2}{2} \sum _{j\in \mathbb{Z}}  \coth \left(\frac{\pi
  (|\omega+2\pi j |)}{2}\right) \text{csch} \left(\frac{\pi  (|\omega+2\pi j |)}{2}\right).
 \] 

The cardinal function, which also decays cubically at infinity (\cite[Theorem 4.3]{RBF}, is then obtained similarly as indicated in the previous example. For practical computing purposes, since the tails of the symbol decay again quickly, we proceed analogously to obtain an approximated version of the cardinal function:
\[\chi_t(x) := \sum _{j\in \mathbb{Z}}  c_j \phi(\left| x-j\right| ), \qquad x\in \mathbb{R}.\]

\end{description}
The plot of the obtained wavelet can be observed in Figure  \ref{canena}.
\subsection{Wavelet based on the Thin-Plate Spline}

Similarly to the previous cases, we now take $\phi(r)= (r^2+ c^2) \log(c^2 + r^2)$ with  $c > 0$. 
The generalised Fourier transform of the function $\phi$ is given by (see \cite{ORT})
\begin{equation}\label{TPSFT}
\hat{\phi}(\omega)=2\sqrt{\pi} \left(\dfrac{2c}{|\omega|}\right)^{\frac32}K_{\frac32}(c|\omega|), \qquad \omega\in \mathbb{R}.
\end{equation}
In fact, this has a particularly simple form because Macdonald functions of integer orders plus half are surprisingly simple, see \cite{GR}:
$$ K_{\frac32}(z)=\frac{\sqrt\pi e^{-z}}{\sqrt{2z}}\Bigl(1+\frac1z\Bigr).$$

\begin{description}

\item[Quasi-Lagrange function]
The expansion of $\hat{\phi}(\omega)$ about the origin is given by (see \cite{HYP})
\[
4\pi \left( \frac{1}{|\omega|^3}-\frac{c^2}{2|\omega|}+\frac{c^3}{3}-\frac{c^4}{8}|\omega|+\cdots   \right).
\]
Therefore, $\hat{\phi}$ presents a singularity of order 3 at the origin.

By considering, for instance, $\mu_k$ as the Fourier coefficients of the trigonometric polynomial $P(\omega)=|\sin \omega |^3$, we obtain the quasi-Lagrange function
\[
\psi(x)=\sum_{k \in \mathbb{Z}} \mu_k \phi (|x-k|)= \sum_{j\in \mathbb{Z}}\frac{12}{\pi
  \left(9- 4j^2\right) \left(1- 4j^2\right)}\phi(|x-2j|).
\] 
The computation of the coefficients is as follows:
\begin{align*}
c_k=&\frac{1}{2\pi} \int_{-\pi}^{\pi} e^{{\rm i}\omega k} |\sin \omega |^3d\omega=\frac{2}{2\pi} \int_{0}^{\pi} \cos(k \omega)\sin^3 \omega  d\omega \\
=& \frac{1}{\pi} \dfrac{\pi \cos(\frac{\pi k}{2})}{2^3 4 B(\frac{5+k}{2}, \frac{5-k}{2})}
=\dfrac{ \Gamma(5)\cos(\frac{\pi k}{2})}{32 \Gamma(\frac{5+k}{2})\Gamma (\frac{5-k}{2})}
= \dfrac{3\cos(\frac{\pi k}{2})}{4 \Gamma(\frac{5+k}{2})\Gamma (\frac{5-k}{2})},
\end{align*}
where we have used \cite{GR} (3.631 8., 8.331 1., 8.334 3., 8.338 1., and 8.384 1.) for the beta function $B$ and for some properties of the Gamma function.

So, the Fourier coefficients are 0 for odd $k.$ For $k= 2m$ we have that
\[ c_{2m}= (-1)^m\dfrac{3}{4\Gamma\Bigl(\frac{5}{2}+ m\Bigr)\Gamma \Bigl(\frac{5}{2}- m\Bigr)}.\]
\[
\Gamma\Bigl(\frac52 + m\Bigr) =\Bigl(\frac32 + m\Bigr)\Bigl(\frac12 + m\Bigr) \Gamma\Bigl(\frac12 + m\Bigr), \qquad
\Gamma\Bigl(\frac52 - m\Bigr) = \Bigl(\frac32 - m\Bigr)(\frac12 - m\Bigr) \Gamma\Bigl(\frac12 - m\Bigr).
\]
Now, by setting $z = 1/2 + m$, we have $1-z = \frac12 - m$ and 
$$\Gamma\Bigl(\frac12 + m\Bigr) \Gamma\Bigl(\frac12 - m\Bigr) = \frac{\pi}{\sin\left(\pi(\frac12 + m)\right)}$$ and we have 
\begin{align*}
c_{2m}=& (-1)^m\dfrac{3\sin\left(\pi(\frac12 + m)\right)}{4\pi (\frac32 + m)(\frac12 + m) (\frac32 - m)(\frac12 - m) }\\
=&\dfrac{3 }{4\pi (\frac94 - m^2)(\frac14 - m^2)}=\dfrac{12}{\pi(9 - 4m^2)(1 - 4m^2)}.
\end{align*}

With this choice, it is obtained that $\psi$ has the asymptotic property
\[
|\psi(x)|={O}((1+|x|)^{-4}),
\]
for large $|x|$, and also that 
\[
\hat{\psi}(\omega)=P(\omega) \hat{\phi}(|\omega|)=1+{O}(|\omega|^2)
\]
for small argument
(see \cite{JJMM} for a similar reasoning). Therefore, the conditions in Theorem \ref{SF} are satisfied with $m=1$ and $\ell=2$.

\item[Cardinal function]
Taking into account (\ref{TPSFT}), the symbol is now given by 
\[
\sigma (\omega)=  2^{5/2} c^{3/2} \sqrt{\pi} \sum_{j\in \mathbb{Z}}   \dfrac{K_{\frac32}(c|\omega+ 2\pi j|)}{|\omega+ 2\pi j|^{3/2}}, \qquad \omega\in\mathbb{R},
\]
which is
\[
\sigma (\omega)=  4\pi \sum_{j\in \mathbb{Z}}   \dfrac{ \exp(-c|\omega+ 2\pi j|)(1+c|\omega+2\pi j|)}{|\omega+ 2\pi j|^{3}}, \qquad \omega\in\mathbb{R},
\]
due to the form of the Macdonald function stated above.
 
 Since the order of the singularity of $\hat{\phi}$ at the origin is 3 and $\hat{\phi}(r)>0$ for $r>0$ (see \cite{HYP}), the asymptotic decay of $\chi (x)$ is identified as 
\[
|\chi (x)|={O}((1+|x|)^{-4}),
\]
for large $|x|$ (\cite[Theorem 4.3]{RBF}). The construction of the cardinal function will follow similar steps as done in the previous examples.


\end{description}

Although perfect from a theoretical point of view, the use of this
function $\phi$ is almost unsuitable for practical purposes. The convergence of the trigonometric polynomial $P(\omega)$'s Fourier series is very slow, which means
that when we truncate it, we will not obtain a good approximation and
therefore errors occur when computing 
the quasi-Lagrange function (see \cite{HYP} to illustrate the use of
RBF on the wrong parity of the dimension).

\subsection{Wavelet based on the Generalized Multiquadric}
We now take $\phi(r)= \sqrt{c^{2d} + r^{2d}}$ with  $c > 0$ and $d\in \mathbb{N}.$ The details of the study of this RBF and its associated quasi-Lagrange functions can been found in \cite{ORT1}. For this example we take $d= 5$  which will produce a very convenient decay of the quasi-Lagrange function for our purposes.

\begin{description}
\item[Quasi-Lagrange function]
This generalised Fourier transform $\hat{\phi}$ presents a singularity of order 6 at the origin since
 \[
 \hat{\phi}(|\omega|) \sim \frac{240}{\omega ^6}, \qquad |\omega|\rightarrow 0.
 \]
In this case the quasi-Lagrange function has the form
\[\psi(x)=\sum_{k=-7}^{7} \mu_k \phi (|x-k| )\]
where 
\[\mu_{-7}=\mu_{7}=-\frac{c^8 q}{110592}+\frac{c^6 p}{82944}+\frac{37}{3628800},
 \quad \mu_{-6}=\mu_{6}=\frac{7 c^8 q}{55296}-\frac{c^6 p}{5184}-\frac{2767}{14515200},\]
\[\mu_{-5}=\mu_{5}=-\frac{91 c^8 q}{110592}+\frac{115 c^6 p}{82944}+\frac{6271}{3628800},
\quad \mu_{-4}=\mu_{4}=\frac{91 c^8 q}{27648}-\frac{31 c^6 p}{5184}-\frac{73811}{7257600},\]
\[\mu_{-3}=\mu_{3}=-\frac{1001 c^8 q}{110592}+\frac{1441 c^6 p}{82944}+\frac{157477}{3628800},
\quad \mu_{-2}=\mu_{2}=\frac{1001 c^8 q}{55296}-\frac{187 c^6 p}{5184}-\frac{1819681}{14515200},\]
\[\mu_{-1}=\mu_{1}= -\frac{1001 c^8 q}{36864}+\frac{1529 c^6 p}{27648}+\frac{286397}{1209600}
\quad \mu_{0}= \frac{143 c^8 q}{4608}-\frac{55 c^6 p}{864}-\frac{353639}{1209600}\]
and \[q=\dfrac{5 \sqrt{5}\pi^2}{2^{\frac{3}{5}}\Gamma(-\frac{8}{5})\left(\Gamma(-\frac{1}{5})\right)^2}; \quad
p = \dfrac{5 \sqrt{5}\pi^2}{2^{\frac{1}{5}} \Gamma(-\frac{6}{5})\left(\Gamma(-\frac{2}{5})\right)^2}.\]
It has been proven in the aforementioned work that the asymptotic behaviour of this quasi-Lagrange function is $|\psi(x)|={O}(|x|^{-11})$ when $|x|$ goes to infinity. The other conditions in Theorem \ref{SF} are satisfied with $m=9$ and $\ell=1$. 

\item[Cardinal function]We will compute the cardinal function directly, similarly to the first example. Specifically, we obtain it from a linear combination of shifts of $\psi$ by imposing the cardinal function to be one at the origin and zero at all other integers (the cardinality property). 

We must mention that, for practical purposes, we have to proceed more carefully here because, unlike in the previous example, the quasi-Lagrange function's coefficient vector does not have compact support in this case. Thus, instead of generating the cardinal function as
\[  
\chi(x)=\sum_{j\in \mathbb{Z}} \lambda_j \psi(x-j),\qquad x\in{\mathbb R},
\] 
we will approximate it by 
\begin{equation}\label{rec} 
\chi_N(x)=\sum_{j\in \mathbb{Z}} \lambda^{\ast}_j \psi_N(x-j),\qquad x\in{\mathbb R}. 
\end{equation}
There is no harm in doing so because of the quick decay of the quasi-Lagrange function $\psi$. These are the steps we follow: First, we define a new function, 
$\psi_N(x)$, by truncating the quasi-Lagrange function $\psi(x)$ on the interval $I=(-N-1, N+1)$, where $N\in\mathbb{N}$. 
Specifically: 
\[ 
\psi_N(x):= 
\begin{cases}
\psi(x) & \: \text{if } x \in I, \\
0&  \: \text{otherwise.}
\end{cases}
\]
We notice that $\psi_N(j)=\psi_N(-j)$ for all $j \in \mathbb{Z}$ due to the symmetry of $\psi.$

Now, taking into account (\ref{rec}) and by imposing the cardinality property $\chi_N(i)=\delta_{i,0}$ for all integers $i \in \mathbb{Z}$, 
we obtain the infinite linear system:
\begin{equation}\label{sys}
\begin{cases} \displaystyle \sum_{j \in \mathbb{Z}} \lambda^{\ast}_j \psi_N(j)= 1, \\[1ex]
\displaystyle \sum_{j \in \mathbb{Z}} \lambda^{\ast}_{i+j} \psi_N(j) = 0, \quad \forall i \in \mathbb{Z}\setminus\{0\}.
\end{cases}
\end{equation}
 In matrix form, system (\ref{sys}) becomes
\begin{equation}\label{sys0}
A {\bf \lambda^{\ast}} = b
\end{equation}
where $b=(\delta_{i,0})_{i \in \mathbb{Z}}$, $A=(a_{ij})_{i,j\in \mathbb{Z}}$ is a symmetric $2N+1$-banded matrix, and 
${\bf \lambda^{\ast}}=(\lambda^{\ast}_i)_{i\in\mathbb{Z}}$ is the vector of unknowns. 
More precisely, $A$ is given by the Toeplitz matrix
\[A = 
\begin{pmatrix} 
 \ddots &\ddots &  \ddots      &  \ddots         & \ddots & \ddots     &  \ddots & \ddots     &\ddots \\
 \ddots        &  \psi_N(N) & \psi_N(N- 1) & \dots &\psi_N(0) & \dots   &\psi_N(N) & 0& \ddots \\
  \ddots    &0     &  \psi_N(N) & \psi_N(N- 1) & \dots &\psi_N(0) & \dots   &\psi_N(N) & \ddots \\
 \ddots  &\ddots &  \ddots      &  \ddots         & \ddots & \ddots     &  \ddots & \ddots     &\ddots
\end{pmatrix}.
\]
We solve this (infinite) system of equations by using the characteristic equation of degree $2N$
\begin{equation}\label{char}
\sum _{j=-N}^{N}\psi_N(|j|) r^{N+j}= 0.
\end{equation}
We point out  that a palindromic polynomial appears in expression (\ref{char}). This means that if $r$ is a root, then so is $1/r$.  Furthermore, since the sum of the coefficients is greater than 0, $r=1$ cannot be a root of  (\ref{char}). Thus, we know that  (\ref{char}) has exactly $N$ roots with modulus less than 1. 
Now, the general solution of the difference equation associated with (\ref{sys}), which is valid for $ |n| \geq N,$ is given by 

\begin{equation}\label{sol}
\lambda_n^{\ast} = \sum _{m=1}^{2N}c_mr_m^{|n|},
\end{equation}
being $r_m\in \mathbb{C}$  the roots of equation (\ref{char}) and $c_m\in \mathbb{C}$ constants to be determined.  Based on what is stated in the previous paragraph, we can write  (\ref{sol}) as
\begin{equation}\label{sol1}
\lambda_n^{\ast} = \sum _{m=1}^{N}c_mr_m^{|n|}+\sum _{m=N+1}^{2N}c_mr_m^{|n|},
\end{equation}
with  $|r_m|<1, i=1, \ldots N.$
Now, we claim that (\ref{sol}) is a solution of (\ref{sys0}) if we can find coefficients $c_m$ verifying
\begin{equation}\label{condi} 
\sum_{j=-N}^{N}\psi_N(|j|) \lambda_{j+k}^{\ast}= \delta_{k,0}, \qquad k=0, \ldots, N-1
\end{equation}
or, equivalently
\[ \sum _{m=1}^{2N} c_m \left( \sum_{j=-N}^{N}\psi_N(|j|) r_m^{|j+k|} \right) = \delta_{k,0}, \qquad k=0, \ldots, N-1.
\]
For $k=0$, the condition given in (\ref{condi}) comes from the first equation of (\ref{sys}). The other $N-1$ conditions (for $k=1, \ldots, N-1$) come from imposing $\lambda_n^{\ast}$ to be a solution for $ |n| < N$.  So, (\ref{condi}) gives us $N$ conditions and we have $2N$ unknowns to solve the linear system. Taking $c_m=0,  i=N+1, \ldots 2N$ we must solve a linear system with $N$ unknowns and $N$ equations. The solution of  (\ref{sol1}), if exists, gives rise  to a  bounded cardinal function.  If we find no solution, we have always the option to change the value of $N$ (the truncation parameter in the interval $I$) in order to get a solution of  (\ref{sol1}). Finally, a good approximation for the cardinal function is given by 
\[ 
\chi_N(x)=\sum _{j \in \mathbb{Z}} \lambda^{\ast}_j \psi_N(x-j),\qquad x\in{\mathbb R}
\]
and  the associated wavelet obtained from this approximation of the cardinal function is 
\[ 
\zeta(x)= \psi(2x)- \chi_N(2x), \qquad x\in \mathbb{R}.
\] 
For $N=9$, the graphic of  the cardinal function and its associated wavelet  are given in Figure \ref{GMd5}.
\begin{figure}[h]
        \centering
        \includegraphics[width=0.44\textwidth]{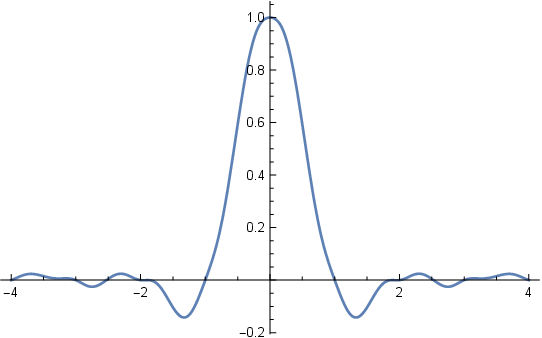}
        \hspace{0.1cm}
          \includegraphics[width=0.44\textwidth]{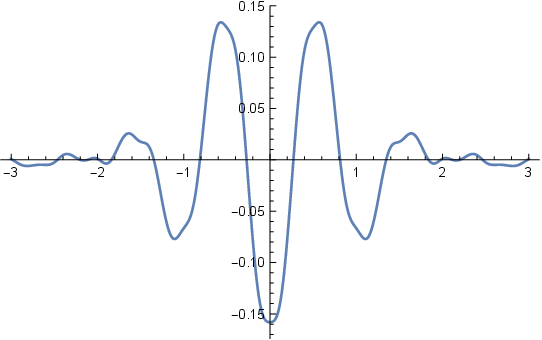}
 \caption{Cardinal function  (left) and wavelet based on $\phi(r)=\sqrt{1 + r^{10}}$ (right).}
 \label{GMd5}  
\end{figure}

\end{description}
\section{Appendix: Fourier transform of $|x| \tanh |x|$}

From \cite[Theorem 5]{HYP} the univariate generalised Fourier transform of the hyperbolic tangent is:
\[
\mathcal{F}\bigl(\tanh x\bigr)(\omega)=-\mathrm{i}\pi\,\text{csch}\,\left(\dfrac{\pi\omega}{2}\right).
\]

By differentiating both terms with respect to $\omega$, we have:
\[
\frac{d}{d\omega}\left[ \int_{-\infty}^{\infty} e^{-\mathrm{i} x \omega} \tanh x\, dx \right]=\frac{d}{d\omega} \left[ -\mathrm{i}\pi\,\text{csch}\,\left(\dfrac{\pi\omega}{2}\right)\right],
\]
that is
\[
-\mathrm{i}  \int_{-\infty}^{\infty} e^{-\mathrm{i} x \omega} x\tanh x \,dx=
 -\mathrm{i} \pi \left(-\dfrac{1}{2}\right) \pi \text{coth}\, \left(\dfrac{\pi\omega}{2}\right) \text{csch}\, \left(\dfrac{\pi\omega}{2}\right).
\]
Therefore:
\[
\mathcal{F}(x \tanh x)(\omega)=-\dfrac{\pi^2}{2} \text{coth}\, \left(\dfrac{\pi\omega}{2}\right) \text{csch}\, \left(\dfrac{\pi\omega}{2}\right).
\]
On the other hand, since $\tanh (-x)=-\tanh x$, it is satisfied that $\vert x \vert \tanh \vert x \vert=x \tanh x$. Therefore we can conclude that:
\[
\mathcal{F}(\vert x \vert \tanh \vert x \vert)(\omega)=-\dfrac{\pi^2}{2} \text{coth}\, \left(\dfrac{\pi\omega}{2}\right) \text{csch}\, \left(\dfrac{\pi\omega}{2}\right).
\]

\affiliationone{
   M.D. Buhmann\\
   Justus-Liebig University, Department of Mathematics, 35392, Giessen,
        Germany
   \email{buhmann@math.uni-giessen.de}}
\affiliationtwo{
   J. J\'odar\\
 University of Jaén, Department of  Mathematics, 23071, Jaén, 
  Spain
   \email{jjodar@ujaen.es}}
\affiliationthree{~} 
\affiliationfour{
   M. L. Rodr\'iguez\\
   University of Granada, Department of Applied Mathematics, 18071, Granada,
   Spain
   \email{miguelrg@ugr.es}}
\end{document}